
\input amssym.def
\input amssym
\input psfig
\magnification=1100
\baselineskip = 0.23truein
\lineskiplimit = 0.01truein
\lineskip = 0.01truein
\vsize = 8.5truein
\voffset = 0.2truein
\parskip = 0.10truein
\parindent = 0.3truein
\settabs 12 \columns
\hsize = 5.8truein
\hoffset = 0.2truein

\setbox\strutbox=\hbox{%
\vrule height .708\baselineskip
depth .292\baselineskip
width 0pt}
\font\caps=cmcsc10

\font\bigtenrm=cmr10 at 14pt

\def\sqr#1#2{{\vcenter{\vbox{\hrule height.#2pt
\hbox{\vrule width.#2pt height#1pt \kern#1pt
\vrule width.#2pt}
\hrule height.#2pt}}}}
\def\square{\mathchoice\sqr46\sqr46\sqr{3.1}6\sqr{2.3}4}

\centerline{\bigtenrm A CHARACTERISATION OF}
\centerline{\bigtenrm LARGE FINITELY PRESENTED GROUPS}

\tenrm
\vskip 14pt
\centerline{MARC LACKENBY}
\vskip 18pt

\centerline{\caps 1. Introduction}
\vskip 6pt

In this paper, we will consider finitely presented groups that have
a finite index subgroup which admits a surjective
homomorphism onto a non-abelian free group.
Gromov called these groups {\sl large} [4].
Large groups have particularly nice properties
(for example, super-exponential subgroup growth).
They also play an important r\^ole in low-dimensional
topology: it is a major conjecture that the
fundamental group of any closed hyperbolic
3-manifold is large.
Our main theorem is a characterisation of
these groups in terms of the existence of a
normal series where successive quotients are
finite abelian groups with sufficiently large
rank and order.

\noindent {\bf Theorem 1.1.} {\sl Let $G$ be a finitely
presented group. Then the following are equivalent:
\item{1.} some finite index subgroup of $G$ admits
a surjective homomorphism onto a non-abelian free group;
\item{2.} there exists a sequence $G_1 \geq G_2 \geq \dots$
of finite index subgroups of $G$, each normal in $G_1$, such that
\itemitem{(i)} $G_i/G_{i+1}$ is abelian for all $i \geq 1$;
\itemitem{(ii)} $\lim_{i \rightarrow \infty} 
((\log [G_i : G_{i+1}]) / [G:G_i]) = \infty$;
\itemitem{(iii)} $\limsup_i (d(G_i/G_{i+1}) / [G:G_i])  > 0$.
}

Here, $d( \ )$ denotes the rank of a group, which
is its minimal number of generators.
Note that condition (i) does not require $G/G_1$
to be an abelian group. Indeed, it cannot, since there
are finitely presented groups that are both perfect
and large. In words, condition (ii) requires 
the order of $G_i/G_{i+1}$ to grow super-exponentially
as a function of $[G:G_i]$. Condition (iii)
asserts that the rank of the
quotients $G_i/G_{i+1}$ grows linearly in the
index $[G:G_i]$. This is the fastest it could
possibly grow, since the Reidermeister-Schreier
theorem implies that the rank of $G_i$ grows
at most linearly in $[G:G_i]$.

The main part of Theorem 1.1, the implication
$(2) \Rightarrow (1)$, is in fact a corollary
of the following stronger, but slightly less elegant, result.
\vfill\eject

\noindent {\bf Theorem 1.2.} {\sl Let $G$ be a
finitely presented group, and suppose that,
for each natural number $i$, there is a
triple $H_i \geq J_i \geq K_i$ of finite index
normal subgroups of $G$ such that
\item{(i)} $H_i/J_i$ is abelian for all $i$;
\item{(ii)} $\lim_{i \rightarrow \infty} 
((\log [H_i : J_i]) / [G:H_i]) = \infty$;
\item{(iii)} $\limsup_i (d(J_i/K_i) / [G:J_i])  > 0$.

\noindent Then $K_i$ admits a surjective homomorphism
onto a free non-abelian group, for infinitely
many $i$.}

This gives $(2) \Rightarrow (1)$ of Theorem 1.1,
via the following argument. We may replace $G$ by
$G_1$, and thereby assume that each $G_i$ is
normal in $G$. If we then set $H_i = G_i$, $J_i = G_{i+1}$
and $K_i = G_{i+2}$, Theorem 1.2 implies that
infinitely many $G_i$ admit a surjective homomorphism
onto a free non-abelian group.

In fact, the statements of Theorems 1.1 and 1.2 are
not the strongest that can be made. The existence
of infinitely many subgroups $G_i$ (in Theorem 1.1)
or infinitely many triples $H_i \geq J_i \geq K_i$
(in Theorem 1.2) is more than one needs to deduce
that $G$ is large. This conclusion still holds
if one replaces hypotheses (ii) and (iii) of Theorem 1.2
by an explicit inequality which relates
$[H_i:J_i]$, $[G:H_i]$, $d(J_i/K_i)$, for some
fixed $i$, and data from a fixed
presentation of $G$. The precise statement of this
result, which is rather unwieldy, appears as Theorem
4.2 in Section 4.

What makes the results of this paper noteworthy is their method of
proof. Despite the fact that these theorems are
purely group-theoretic, their proof uses very little
algebra. Instead, the geometry and topology
of finite Cayley graphs play a central r\^ole.
This is the second in a pair of papers that
exploit these type of arguments. The first [5] related
Property $(\tau)$, the rank of finite index subgroups and
their possible decomposition as a graph of groups.
Using this relationship, the proof of a weaker form of
Theorem 1.1 was sketched.

Large groups were studied
by Baumslag and Pride [1] who showed that groups
with a presentation having (at least) two more
generators than relations are of this form.
This is an easy consequence of Theorem 1.1.

\noindent {\bf Corollary 1.3.} {\sl Let $G$
be a group having a presentation with at least
two more generators than relations. Then $G$
has a finite index subgroup that admits a
surjective homomorphism onto a free non-abelian
group.}

Baumslag and Pride also conjectured [2] that
when a group $G$ has a presentation with one more
generator than relation, but one of the relations
is a proper power, then $G$ is large. This was proved 
independently by Gromov [4] using bounded cohomology 
and by St\"ohr [10] using a direct algebraic argument. 
Again, this is a straightforward consequence of Theorem 1.1.

\noindent {\bf Corollary 1.4.} {\sl Let
$G$ be a group having a presentation with
one more generator than relation. Suppose
that one of these relations is proper power.
Then $G$ has a finite index subgroup that
admits a surjective homomorphism onto a free non-abelian
group.}

We will prove these corollaries in \S5.

We will also show that conditions (i) and (ii) of Theorem 1.1 are,
in fact, equivalent to the existence of
a finite index subgroup with infinite
abelianisation.

\noindent {\bf Theorem 1.5.}
{\sl Let $G$ be a finitely
presented group. Then the following are equivalent:
\item{1.} some finite index subgroup of $G$ 
has infinite abelianisation;
\item{2.} there exists a sequence $G_1 \geq G_2 \geq \dots$
of finite index subgroups of $G$, each normal in $G_1$, such that
\itemitem{(i)} $G_i/G_{i+1}$ is abelian for all $i \geq 1$;
\itemitem{(ii)} $\lim_{i \rightarrow \infty} 
((\log [G_i : G_{i+1}]) / [G:G_i]) = \infty$.
}

The proof of this result, which is topological and
rather straightforward, is given in \S6.

I am grateful to Peter Shalen for a discussion
at the early stages of this project, when he
suggested that a result along the lines of Theorem 1.1
might be true. I am also grateful to Alex
Lubotzky for some useful conversations about
Property $(\tau)$, when we discussed conditions
similar to those in Theorem 1.1.

\vskip 18pt
\centerline{\caps 2. The forwards implication}
\vskip 6pt

In this section, we will prove $(1) \Rightarrow (2)$
of Theorem 1.1.
Let $G_1$ be the finite index subgroup of $G$ that admits a surjective
homomorphism onto a non-abelian free group $F$. Define
recursively the following finite index normal
subgroups of $F$. Set $L_1 = F$, and, for $i \geq 1$,
let $L_{i+1} = [L_i,L_i] (L_i)^i$. In other words,
$L_{i+1}$ is the group generated by the commutators
and the $i^{\rm th}$ powers of $L_i$. Thus, $\{ L_i \}$
is a nested sequence of subgroups of $F$. Each is
characteristic in the preceding one, and hence
each is characteristic in $F$, and is therefore
normal in $F$. Each $L_i$ is free with rank
$d(L_i) = (d(F) - 1)[F:L_i] + 1$, and $L_i/L_{i+1}$
is isomorphic to $({\Bbb Z}/i{\Bbb Z})^{d(L_i)}$. 

Set $G_i$ to be the inverse image of $L_i$ in $G$,
which is a normal subgroup of $G_1$. Then
$G_1/G_i$ is isomorphic to $F/L_i$, and $G_i/G_{i+1}$
is isomorphic to $L_i/L_{i+1}$,
which is abelian, verifying (i) of Theorem 1.1.
To check (iii), note that
$$d(G_i / G_{i+1})= d(L_i / L_{i+1}) = d(L_i) = 
(d(F) - 1)[F:L_i] + 1
= (d(F) - 1)[G_1:G_i] + 1.$$
Similarly, to establish (ii), we note that
$$\log [G_i:G_{i+1}] = \log ( i^{d(L_i)} )
= d(L_i) \log i > [G_1:G_i] \log i.$$

\vskip 18pt
\centerline{\caps 3. The width and Cheeger constant of finite graphs}
\vskip 6pt

Many of the ideas behind this paper arise from
the theory of Property $(\tau)$. This is a particularly
useful group-theoretic concept, introduced by Lubotzky and Zimmer [8],
that can be defined using graph
theory, representation theory or differential
geometry. We concentrate on the former approach.
The {\sl Cheeger constant} of a finite graph $X$, denoted
$h(X)$, is defined to be
$$\min \left\{ {|\partial A| \over |A|}:
A \subset V(X) \hbox{ and } 0 < |A| \leq {|V(X)| \over 2} \right\}.$$
Here, $V(X)$ denotes the vertex set of $X$, and,
for a subset $A$ of $V(X)$, $\partial A$ denotes
the set of edges with one endpoint in $A$ and one
not in $A$. Informally, having small Cheeger constant
is equivalent to the existence of a `bottleneck'
in the graph. (See Figure 1.)

\vskip 18pt
\centerline{\psfig{figure=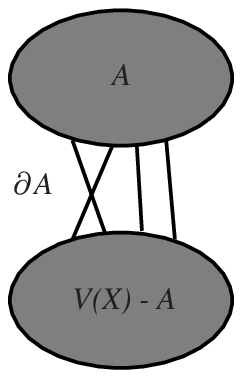}}
\vskip 18pt
\centerline{Figure 1.}

Let $G$ be a group with a finite generating set $S$.
Let $\{ G_i \}$ be a collection of finite index
normal subgroups. We denote the Cayley graph of
$G/G_i$ with respect to $S$ by $X(G/G_i;S)$.
The group $G$ is said to have {\sl Property $(\tau)$}
with respect to $\{ G_i \}$ if the Cheeger constants
$h(X(G/G_i;S))$ are bounded away from zero.
This property turns out not to depend on the
choice of finite generating set $S$.

Whether or not a given group and a collection
of finite index subgroups have Property $(\tau)$
is a subtle and often difficult question.
The following theorem of Lubotzky and Weiss [7]
gives a necessary condition for a group to 
have Property $(\tau)$. This is not, in fact,
how Lubotzky and Weiss stated their result
(which appears as Theorem 3.6 of [7]),
but this formulation can readily be deduced
from their argument.

\noindent {\bf Theorem 3.1.} {\sl Suppose that a finitely
generated group $G$ has Property $(\tau)$ with respect
to a collection ${\cal C}$ of finite index normal subgroups.
Then there is a constant $c$ with the following
property. If $J$ is a member of ${\cal C}$,
and $J$ is contained in a normal subgroup $H$ of $G$
such that $H/J$ is abelian, then $|H/J| < c^{[G:H]}$.}

Thus, conditions (i) and (ii) of Theorem 1.2
imply that $G$ does not have
Property $(\tau)$ with respect to $\{ J_i \}$.
We will actually need to establish a stronger
version of Theorem 3.1. Instead of relating to
the Cheeger constant of $X(G/J_i;S)$ (for
some finite generating set $S$), we need to
consider a related geometric invariant
of $X(G/J_i;S)$, its width, which is defined
as follows.

Let $X$ be a finite graph. Consider a linear ordering
on its vertices. For
$0 \leq n \leq |V(X)|$, let $D_n$ be the
first $n$ vertices. The
{\sl width} of the ordering is defined to
be $\max_n |\partial D_n|$. The {\sl
width} of the graph is the minimal width of
any of its orderings, and is denoted $w(X)$.

This notion is inspired by a useful concept from
the theory of knots and 3-manifolds, known as thin position [9],
which was first introduced by Gabai [3]. We now
develop this analogy (which is not essential
for an understanding of the remainder
of the paper). One may imagine the graph
$X$ embedded in ${\Bbb R}^3$, with its
vertices all at distinct heights, and with
its edges realised as straight lines. The height
of the vertices specifies a linear ordering on them.
The width of this ordering can be interpreted
geometrically, as follows. Imagine a 1-parameter family of
horizontal planes, parametrised by their heights
which increase monotonically from $-\infty$ to $\infty$.
The width of the ordering is equal to the
maximal number of intersections between the graph
and any of these planes. (See Figure 2 for
an example.) Thus, to determine the width of
$X$, one should aim to find the most efficient
embedding in ${\Bbb R}^3$: the one that minimises
the width of the associated ordering. This is
highly analogous to thin position for knots in
${\Bbb R}^3$, where one aims to isotope the
knot until a similar notion of width is minimised.

\vskip 18pt
\centerline{\psfig{figure=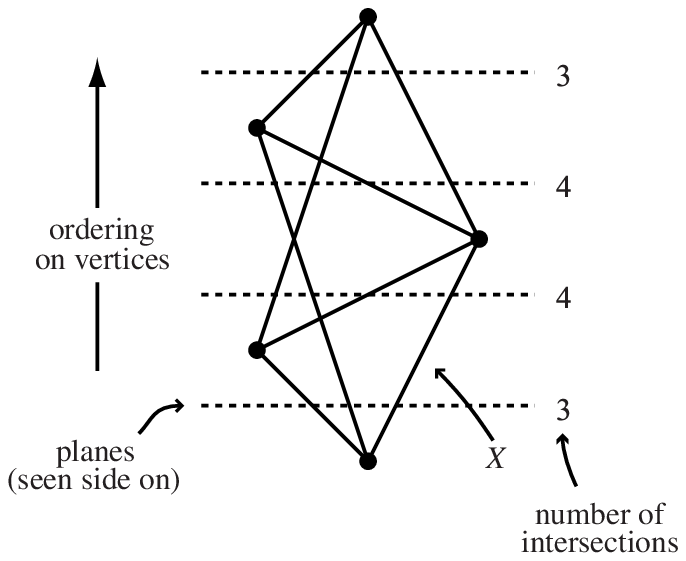}}
\vskip 18pt
\centerline{Figure 2.}

There is a relationship between the width of
a graph and its Cheeger constant.
In an ordering on $V(X)$ of minimal width,
consider $D_n$ where $n = \lfloor |V(X)|/2 \rfloor$.
By computing $|\partial D_n|/|D_n|$, we deduce that
$$h(X) \leq {w(X) \over \lfloor |V(X)|/2 \rfloor}.$$
Hence, the following gives a stronger
version of Theorem 3.1. It gives a condition
guaranteeing that certain Cayley graphs have
width which is asymptotically smaller than
their number of vertices.

\noindent {\bf Theorem 3.2.} {\sl Let $G$ be 
a group with a finite generating set $S$. Suppose that,
for each natural number $i$, there is a pair
$H_i \geq J_i$ of finite index normal
subgroups of $G$, such that
\item{(i)} $H_i/J_i$ is abelian for all $i$;
\item{(ii)} $\lim_{i \rightarrow \infty} 
((\log [H_i : J_i]) / [G:H_i]) = \infty$.

\noindent Then $w(X(G/J_i;S))/[G:J_i] \rightarrow 0$.
}

Note that (i) and (ii) are precisely those in Theorem 1.2.

The proof we give of this result follows the argument of
Lubotzky and Weiss in their proof of Theorem 3.1.
The following lemma allows us to estimate the width of
$X(G/J_i;S)$ in terms of the width
of a Cayley graph of $H_i / J_i$. This will be
useful, since $H_i / J_i$ is abelian, and,
later, we will analyse the width of Cayley
graphs of finite abelian groups.

\vfill\eject

\noindent {\bf Lemma 3.3.} {\sl Let $G$ be a group with
a finite generating set $S$, and let $H_i \geq J_i$
be finite index normal subgroups of $G$. Let $\Sigma$ be a
generating set for $H_i$ arising from $S$ by the
Reidermeister-Schreier process. Then
$$w(X(G/J_i;S)) \leq w(X(H_i/J_i;\Sigma))+2 |S| [G:H_i].$$}

\noindent {\sl Proof.} We first recall the Reidermeister-Schreier
process. Pick a presentation for $G$ with generating
set $S$, but possibly an infinite number of relations.
(We are not assuming here that $G$ is finitely
presented.)
Build the associated 2-complex $C$, by starting with
a bouquet of $|S|$ circles and attaching on a
2-cell for each relation. Then $\pi_1(C)$ is
isomorphic to $G$. Let $C' \rightarrow C$
be the covering corresponding to the subgroup $H_i$, 
so that $\pi_1(C')$ is isomorphic to
$H_i$. Pick a maximal tree $T$ in the 1-skeleton
of $C'$. Collapsing this tree to a point gives
a new 2-complex $\overline C'$. Its 1-cells give the required generating
set $\Sigma$ for $H_i$.

Let $C'' \rightarrow C'$ and $\overline C'' \rightarrow \overline C'$
be the coverings corresponding to $J_i$.
The inverse image of $T$ in $C''$ is
a forest $F$. If one were to collapse each component
of this forest to a single vertex, we would
obtain $\overline C''$. The 1-skeletons of $C''$
and $\overline C''$ are, respectively, 
$X(G/J_i;S)$ and $X(H_i/J_i;\Sigma)$.
(See Figure 3.)

\vskip 18pt
\centerline{\psfig{figure=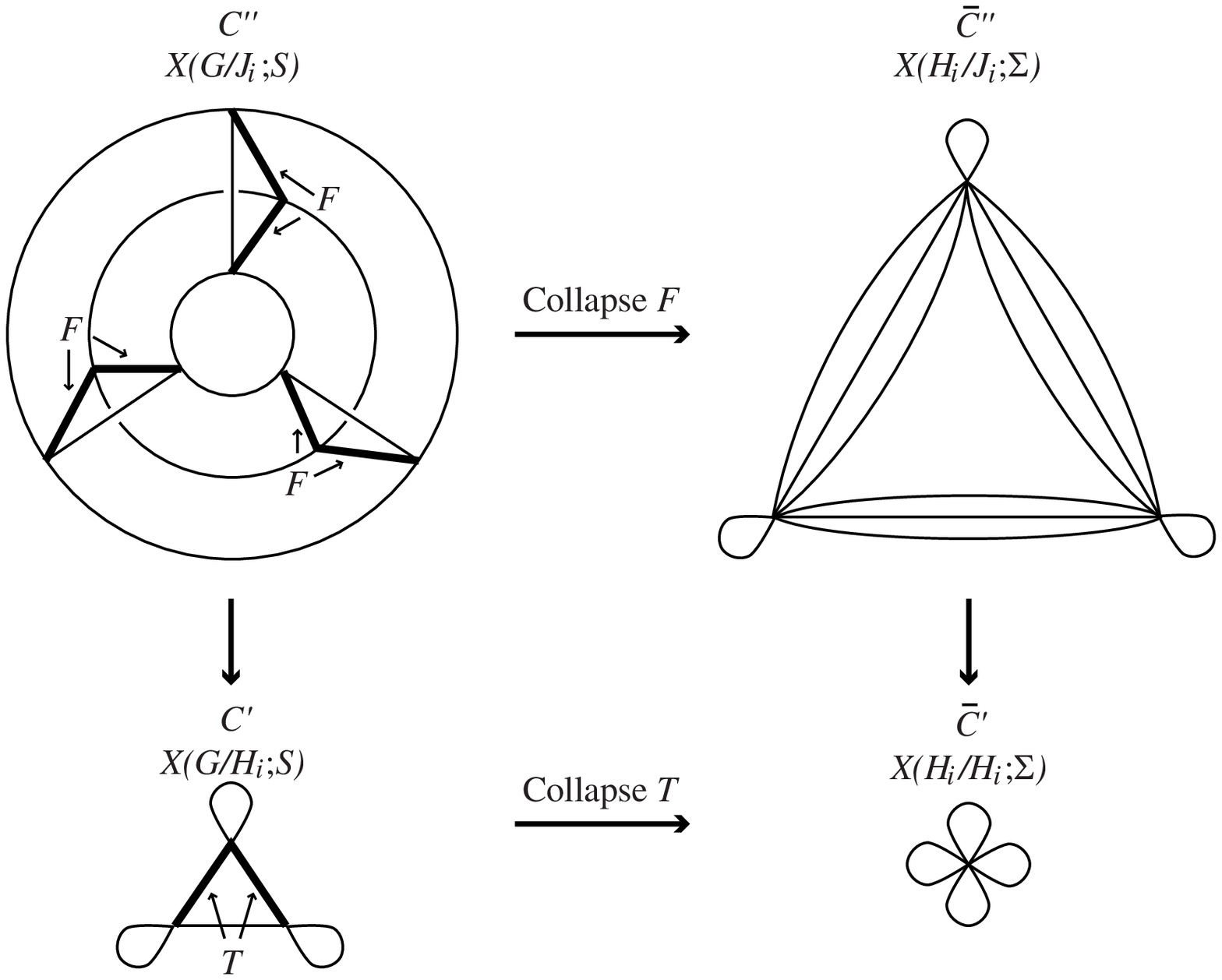,width=4in}}
\vskip 18pt
\centerline{Figure 3.}

Consider a minimal width ordering of $V(\overline C'')$.
From this, one can construct an ordering
of $V(C'')$, as follows. The ordering on $V(\overline C'')$
specifies an ordering on the components of $F$.
We therefore order $V(C'')$ by ordering
the vertices of its first tree in some
way, then the vertices of its second tree,
and so on.

For any integer $n$, with $0 \leq n \leq |V(C'')|$,
let $D_n$ be the first $n$ vertices of $C''$.
Suppose that the $n^{\rm th}$ vertex of
$C''$ lies in a component $\tilde T$ of $F$.
Then, the edges in $\partial D_n$ consist
of edges joining trees other than $\tilde T$, and edges
with at least one endpoint in $\tilde T$. There are at most 
$w(X(H_i/J_i;\Sigma))$ edges of the first type
and at most $2 |S| [G:H_i]$ edges of the second type.
Thus the width of $C''$ is at most
the sum of these quantities. $\square$

As a result of the above lemma, we concentrate
on the Cayley graph of the abelian group
$H_i/J_i$. The following lemma provides a useful
upper bound on its width.

\noindent {\bf Lemma 3.4.} {\sl Let $A$ be a finite
abelian group with finite generating set $\Sigma$.
Then 
$$w(X(A; \Sigma)) \leq 
{6 |\Sigma| |A| \over \lfloor (|A|-1)^{1/|\Sigma| }\rfloor}.$$}

\noindent {\sl Proof.} We will construct an efficient
ordering of the vertices of $X(A; \Sigma)$ by placing
them on the unit circle in a suitable way.
Give the circle a group structure, by identifying
it with ${\Bbb C}^\times$, the 
multiplicative group of complex numbers with modulus one.
Any homomorphism $\phi \colon A \rightarrow {\Bbb C}^\times$
determines a point in $({\Bbb C}^\times)^{|\Sigma|}$,
given by the $|\Sigma|$-tuple $(\phi(s): s \in \Sigma)$.
Define $c$ to be $\lfloor (|A|-1)^{1/|\Sigma|} \rfloor$,
which is the denominator in the upper bound on
$w(X(A; \Sigma))$ that we are trying to establish. 
We may assume that $c$ is a positive integer; otherwise,
there is nothing to prove. Divide
the circle ${\Bbb C}^\times$ into $c$
equal arcs. This determines a decomposition
of $({\Bbb C}^{\times})^{|\Sigma|}$ into
$c^{|\Sigma|} < |A|$ boxes. There are precisely $|A|$ distinct
homomorphisms $A \rightarrow {\Bbb C}^\times$,
and hence two distinct homomorphisms
are sent to the same box. Their
quotient is a non-trivial homomorphism
$\phi \colon A \rightarrow {\Bbb C}^\times$,
such that $|{\rm arg}(\phi(s))| \leq 2\pi/c$
for all $s \in \Sigma$. Here, we are taking
arguments to lie in the range $(-\pi, \pi]$.

Let $\sigma$ be the element of $A / {\rm Ker}(\phi)$,
such that $\phi(\sigma)$ has smallest positive
argument. Then $\sigma$ is a generator for
$A / {\rm Ker}(\phi)$. Let $N$ be its order.
Note that, for any $s \in \Sigma$ that does not lie in ${\rm Ker}(\phi)$,
$N \geq 2 \pi / |{\rm arg}(\phi(s))|
\geq c$. Since this is true for at least
one $s \in \Sigma$, we deduce that $N \geq c$.

We place an order on the vertices of $A$,
as follows. First order the vertices in ${\rm Ker}(\phi)$ in
some way, then order the coset $\phi^{-1}(\sigma)$,
then $\phi^{-1}(\sigma^2)$, and so on. (See Figure 4.)

\vskip 18pt
\centerline{\psfig{figure=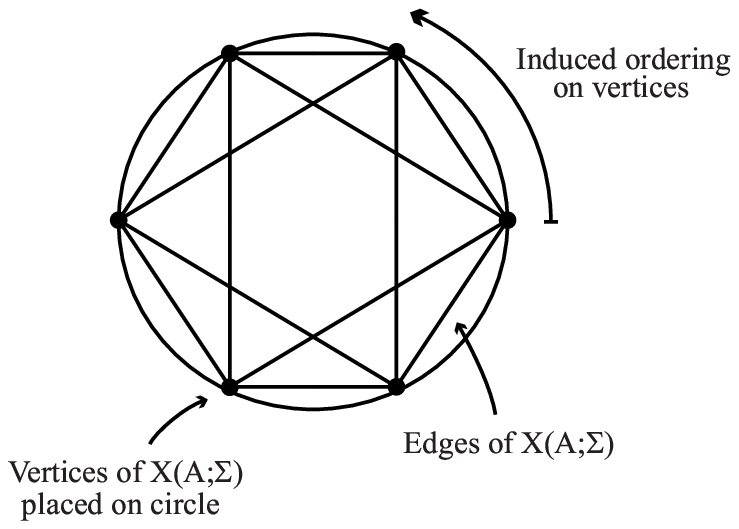}}
\vskip 18pt
\centerline{Figure 4.}

Note that any edge in $X(A; \Sigma)$ joins
vertices who images in $A/{\rm Ker}(\phi)$
differ by $\sigma^i$, where $|i|$ is at most $N/c$. 
Let $D_n$ be a subset of the vertices of
$X(A; \Sigma)$, as in
the definition of width. If we consider
the initial vertex of an edge in $\partial D_n$,
its image in $A/{\rm Ker}(\phi)$ can
take at most $4(N/c)+2$ possible values.
This is at most $6(N/c)$, since $N/c \geq 1$.
Hence, there are at most $6(N/c) |{\rm Ker}(\phi)|
= 6|A|/c$ such vertices. Since $|\Sigma|$
oriented edges emanate from any vertex,
there are
therefore at most $6 |\Sigma| |A| / c$ edges
in $\partial D_n$. This gives the required upper
bound on the width of $X(A; \Sigma)$.
$\square$

\noindent {\sl Proof of Theorem 3.2.}
By Lemma 3.3,
$${w(X(G/J_i;S)) \over [G:J_i]}
\leq {w(X(H_i/J_i;\Sigma)) \over [G:J_i]}
+ {2 |S| [G:H_i] \over [G:J_i]}.$$
By assumption (ii), $[H_i: J_i]$ grows super-exponentially
as a function of $[G:H_i]$. In particular,
$[H_i:J_i]$ tends to infinity,
and so the second term tends to zero.
For the first term, assumption (i) and Lemma 3.4 imply that
$${w(X(H_i/J_i;\Sigma)) \over [G:J_i]}
\leq {6 |\Sigma| \over
[G:H_i] \lfloor  ([H_i:J_i] -1)^{1/|\Sigma|} \rfloor}.$$
From the Reidermeister-Schreier process,
$|\Sigma|$ grows linearly in $[G:H_i]$.
So, $|\Sigma|/[G:H_i]$ is bounded
above and, by assumption (ii),
$([H_i:J_i]-1)^{1/|\Sigma|}$ tends to infinity.
$\square$

\vskip 18pt
\centerline{\caps 4. Proof of the main theorem}
\vskip 6pt

Theorem 1.2 follows immediately from Theorem 3.2
and the following result.

\noindent {\bf Theorem 4.1.} {\sl Let $G$ be
a finitely presented group with a finite
generating set $S$. Suppose that, for
each natural number $i$, there is a
pair $J_i \geq K_i$ of finite index
normal subgroups of $G$, such that
\item{(I)} $w(X(G/J_i;S))/[G:J_i] \rightarrow 0$;
\item{(II)} $\limsup_i (d(J_i/K_i) / [G:J_i])  > 0$.

\noindent Then, for infinitely many $i$,
$K_i$ admits a surjective homomorphism onto a free non-abelian
group.}

Note that in the above theorem, we no longer 
need the hypothesis that any quotient groups
are abelian. 

The proof we give of Theorem 4.1 resembles some of the arguments in [5].
Suppose that (I) and (II) of the theorem hold.
Let $C$ be a finite 2-complex having 
fundamental group $G$, arising from a finite
presentation of $G$ with generating set $S$. Thus, $C$ has
a single vertex, and an oriented edge for
each element of $S$. Let $L$ be
the sum of the lengths of the relations in this
presentation, which we may assume is at least one.

Let $C_i \rightarrow C$ be the covering 
corresponding to $J_i$, and let $X_i = X(G/J_i; S)$ be the
1-skeleton of $C_i$. 
Set $\ell$ to be 
$$\limsup_i {d(J_i/K_i) \over [G:J_i]},$$
which by (II) is positive. For infinitely many $i$,
$${d(J_i/K_i) \over [G:J_i]} > {3 \ell \over 4}.\eqno{(1)}$$
By (I),
$w(X_i)/[G:J_i] \rightarrow 0$.
Hence, provided $i$ is sufficiently large,
$${w(X_i) \over [G:J_i]} < {\ell \over 8 L}.\eqno{(2)}$$
Also, for all large $i$,
$$[G:J_i] \geq {8 (2 |S| + L^2) \over \ell}.\eqno{(3)}$$
We now fix some $i$ so that the inequalities 
(1), (2) and (3) hold. 

We will show that
there is a surjective homomorphism from $K_i$ onto
a free non-abelian group. The details of
this argument are a little complicated, but
the main idea is fairly simple. We use a
minimal width ordering on the vertices of $X_i$
to divide $C_i$ into two pieces $A$ and $B$
that are, roughly, equally `big'. More precisely,
the images of their fundamental groups in $J_i/K_i$
both have large rank. (In fact, $A$ and $B$ might
not be connected, and so it might not make sense
to refer to their fundamental group, but we will
ignore this point for the moment.) Since the
width of $X_i$ is small, the fundamental group
of $A \cap B$ has small rank, as therefore
does its image in $J_i/K_i$. So, when the
decomposition of $C_i$ into $A$ and $B$ is
lifted to the covering space corresponding to
$K_i$, the lifted decomposition is modelled on
a graph with negative Euler characteristic. The
desired conclusion, that $K_i$ surjects onto a
non-abelian free group, follows immediately.

We now give this argument in detail.
Pick a minimal width ordering on the vertices
of $X_i$. We will use this to define
a height function $f \colon C_i \rightarrow {\Bbb R}$.
On the vertices, simply let $f$ agree with
the ordering function. Extend linearly
over the 1-cells of $C_i$. Then we extend
$f$ continuously over the interior of
each 2-cell, so that it is a Morse function
there, with the following properties.
(There is one exceptional case: when all the vertices in
the boundary of a 2-cell are in fact the
same 0-cell in $C_i$ and hence have the same
height, then we define $f$ to be
constant on that 2-cell.)
We can ensure that $f$ has no maxima or minima
in the interior of any 2-cell, and that
its critical values avoid $n + {1 \over 2}$,
for each integer $n$. We can also ensure
that, whenever a 2-cell does not contains the vertex
of height $n$ in its boundary, then that
2-cell contains no critical points with values between
$n-{1 \over 2}$ and $n + {1 \over 2}$.

Let $n$ be some integer between $0$ and $|V(X_i)|$.
We will focus on a single value of $n$ later.
We can now use $n$ and $f$ to decompose $C_i$ into
two subsets $A_n = f^{-1}(-\infty, n+ {1 \over 2}]$
and $B_n = f^{-1}[n+{1 \over 2}, \infty)$.
Let $D_n$ be the vertices of $C_i$ that
lie in $A_n$; these are precisely the
first $n$ vertices of the ordering.

Note that $A_n \cap B_n$ is a 1-complex,
with precisely one 0-cell in the interior
of each edge of $\partial D_n$, and no
other 0-cells. Hence, there are exactly
$|\partial D_n|$ 0-cells of $A_n \cap B_n$.
This therefore is an upper bound on the
number of components of $A_n \cap B_n$, and
hence on the number of components of $A_n$ and
on the number of components of $B_n$ (provided
$n \not= 0, |V(X_i)|$).
Note that $|\partial D_n| \leq w(X_i) < 
{{1 \over 8}} \ell [G:J_i]$, by inequality (2).

The 1-cells of $A_n \cap B_n$ are arcs properly 
embedded in the 2-cells of $C_i$.
We claim that $A_n \cap B_n$ has at most $|\partial D_n| L /2$
1-cells. Note that the total number of
times the 2-cells of $C_i$ run over any 
1-cell of $C_i$ is at most $L$. This is an
upper bound for the valence of each 0-cell
of $A_n \cap B_n$. There are precisely
$|\partial D_n|$ 0-cells of $A_n \cap B_n$,
and hence we obtain the required bound on
the number of 1-cells of $A_n \cap B_n$.
Note that this bound $|\partial D_n| L /2$
is at most $w(X_i) L /2$, which is less than
${{1 \over 8}} \ell [G:J_i]$, by inequality (2).

We construct a graph $Y_n$, as follows. Each
vertex corresponds to a component of $A_n$
or $B_n$, and is labelled $A_n$ or $B_n$ as appropriate. 
Each edge corresponds to a component
of $A_n \cap B_n$. Incidence in the graph is
defined by topological incidence.
Since $A_n \cap B_n$ has a small regular
neighbourhood that is homeomorphic to $(A_n \cap B_n) \times I$,
we may define a collapsing map $C_i \rightarrow Y_n$
that sends each component of $C_i - ((A_n \cap B_n) \times I)$
to the corresponding vertex of $Y_n$, and that
sends each component of $(A_n \cap B_n) \times I$
to the relevant edge of $Y_n$, via projection
onto the second factor of the product.
However, $Y_n$ need not be very interesting:
it may only be a single edge, for example.
But a similar construction in the covering
space of $C_i$ corresponding to the subgroup
$K_i$ will induce the required surjective homomorphism
from $K_i$ onto a free non-abelian group.

\vskip 18pt
\centerline{\psfig{figure=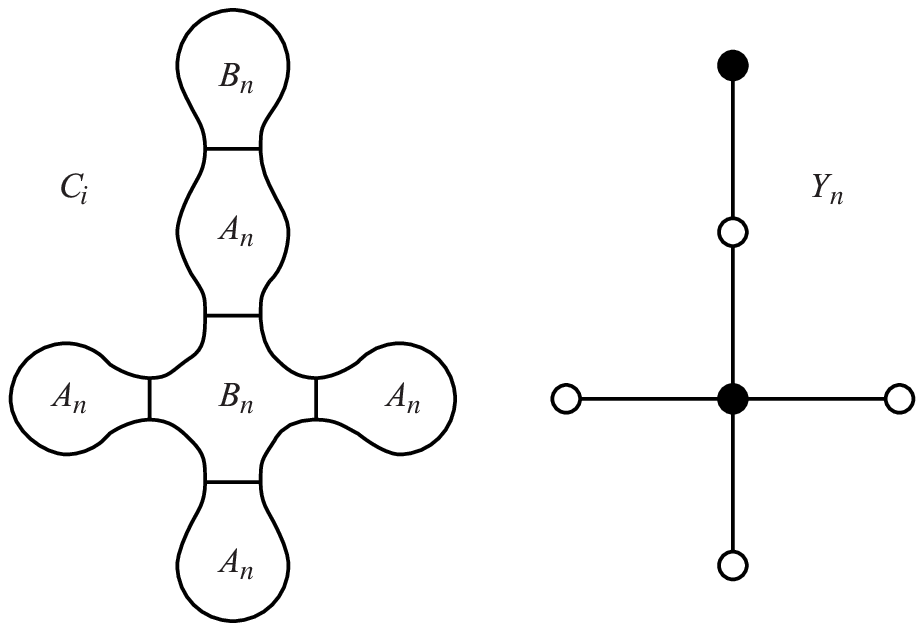}}
\vskip 18pt
\centerline{Figure 5.}

Define the {\sl weight} $wt( \ )$ of
a vertex or edge, corresponding to a
component $E$ of $A_n$, $B_n$ or $A_n \cap B_n$,
to be $d(i_\ast \pi_1 E / i_\ast \pi_1 E \cap K_i)$,
where $i_\ast \colon \pi_1E \rightarrow \pi_1 C_i$ is
the homomorphism induced by inclusion. 
(Note that the choice of basepoints
for $\pi_1E$ does not affect this quantity.)
Define the {\sl weight}
$wt( \ )$ of $A_n$, $B_n$ or $A_n \cap B_n$ to be
the sum of the weights of its components.
Note that $wt(A_n \cap B_n)$
is at most the number of 1-cells of
$A_n \cap B_n$, which we have already
established is less than ${{1 \over 8}} \ell [G:J_i]$.

\noindent {\sl Claim.} There is some $n$
such that $wt(A_n)$ and $wt(B_n)$ are
each at least ${1 \over 4} \ell [G:J_i]$.

Define
$$\eqalign{
Q &= \{ n: wt(A_n) < {\textstyle{1 \over 4}} \ell [G:J_i] \} \cr
R &= \{ n: wt(B_n) < {\textstyle{1 \over 4}} \ell [G:J_i] \}. \cr}$$
The aim is to show that $Q$ and $R$ do not
cover the interval from $0$ to $|V(X_i)|$.
Since $J_i$ is generated by the elements of $i_\ast \pi_1E$,
where $E$ runs over the components of $A_n$ and $B_n$, 
together with at most one generator
for each component of $A_n \cap B_n$,
we have the inequality
$$d(J_i/K_i) \leq
wt(A_n) + wt(B_n) + |A_n \cap B_n|,$$
and so, by the above bound on $|A_n \cap B_n|$ and inequality (1),
$$wt(A_n) + wt(B_n) \geq d(J_i/K_i) - 
{\textstyle{1 \over 8}} \ell [G:J_i] > 
{\textstyle{5 \over 8}} \ell [G:J_i].$$
Hence, when $n \in Q$, $wt(B_n) > {\textstyle{3 \over 8}} \ell [G:J_i]$ 
and, when $n \in R$, $wt(A_n) > {\textstyle{3 \over 8}} \ell [G:J_i]$. 
Thus, $Q$ and $R$ are disjoint.
Note that $0 \in Q$ and $|V(X_i)| \in R$.
Hence, the only way $Q$ and $R$ could cover
the interval from $0$ to $|V(X_i)|$ is
if $n \in Q$ and $n+1 \in R$ for some $n$.
This implies that
$$wt(A_n) < {\textstyle{1 \over 4}}\ell [G:J_i]
\qquad wt(A_{n+1}) > {\textstyle{3 \over 8}}\ell[G:J_i].$$
However, we shall now show that $wt(A_{n})$ and $wt(A_{n+1})$
differ by at most $2|S| + L^2$. Since
$2|S| + L^2 \leq {\textstyle{1 \over 8}}\ell [G:J_i]$,
by inequality (3), this will provide a contradiction.
It is clear that $A_{n+1}$ contains $A_n$. 
Only one vertex of $C_i$ lies in $A_{n+1}$ but
not $A_n$; this is the $(n+1)^{\rm st}$ vertex of the
ordering, called $x$, say. In
the 1-cells and 2-cells of $C_i$ that are disjoint
from $x$, $A_n$ and $A_{n+1}$ differ only by a small
collar. Hence, we need only focus on the 1-cells
and 2-cells that are adjacent to $x$. There are
at most $2 |S|$ of these 1-cells, and at most $L$
2-cells. In each 2-cell adjacent to $x$, we can obtain $A_{n+1}$
from $A_n$ by adding on a collection of discs
that intersect $A_n$ and the boundary of the 2-cell in
a total of at most $L$ arcs. Hence, it is clear
that the weights of $A_{n+1}$ and $A_n$ differ
by at most $2|S| + L^2$. This proves the claim.

We now fix $n$ as in the claim, and abbreviate
$A_n$, $B_n$, $A_n \cap B_n$ and $Y_n$ to $A$, $B$, $A \cap B$ and $Y$.
For any vertex $u$ of $Y$, let $g(u)$ denote its
weight minus the total weight of the edges
to which it is incident. Then the sum of $g(u)$,
over all vertices $u$ of $Y$ labelled $A$,
is $wt(A) - wt(A \cap B)$, which is more than
${{1 \over 8}} \ell [G:J_i]$, and this is more
than the number of vertices labelled $A$.
Hence, there is some vertex $u$ labelled $A$, with
$g(u) > 1$. Similarly, there is some vertex
$v$ labelled $B$ with $g(v) > 1$.
Let $P$ be an embedded path in $Y$ from $u$ to $v$.

Let $p \colon \tilde C_i \rightarrow C_i$ be the covering
corresponding to the subgroup $K_i$. The decomposition of $C_i$ into $A$
and $B$ pulls back to form a similar decomposition
of $\tilde C_i$. We obtain a similar graph
$\tilde Y$. The covering map $p \colon \tilde C_i
\rightarrow C_i$ induces a map of graphs
$\tilde Y \rightarrow Y$. (See Figure 6.) Let $\tilde P$
be the inverse image of $P$ in $\tilde Y$.

The valence of each vertex of $\tilde P$ is at
least that of its image in $P$. When this
image is not an endpoint of $P$, the
valence is therefore at least two.
We shall show that each inverse image of
$u$ and $v$ has at least three edges of 
$\tilde P$ emanating from it. Let $e$ be the edge of
$P$ incident to $u$. Let $U$ and $E$ be
the components of $A$ and $A \cap B$
corresponding to $u$ and $e$. 

\vfill\eject
\centerline{\psfig{figure=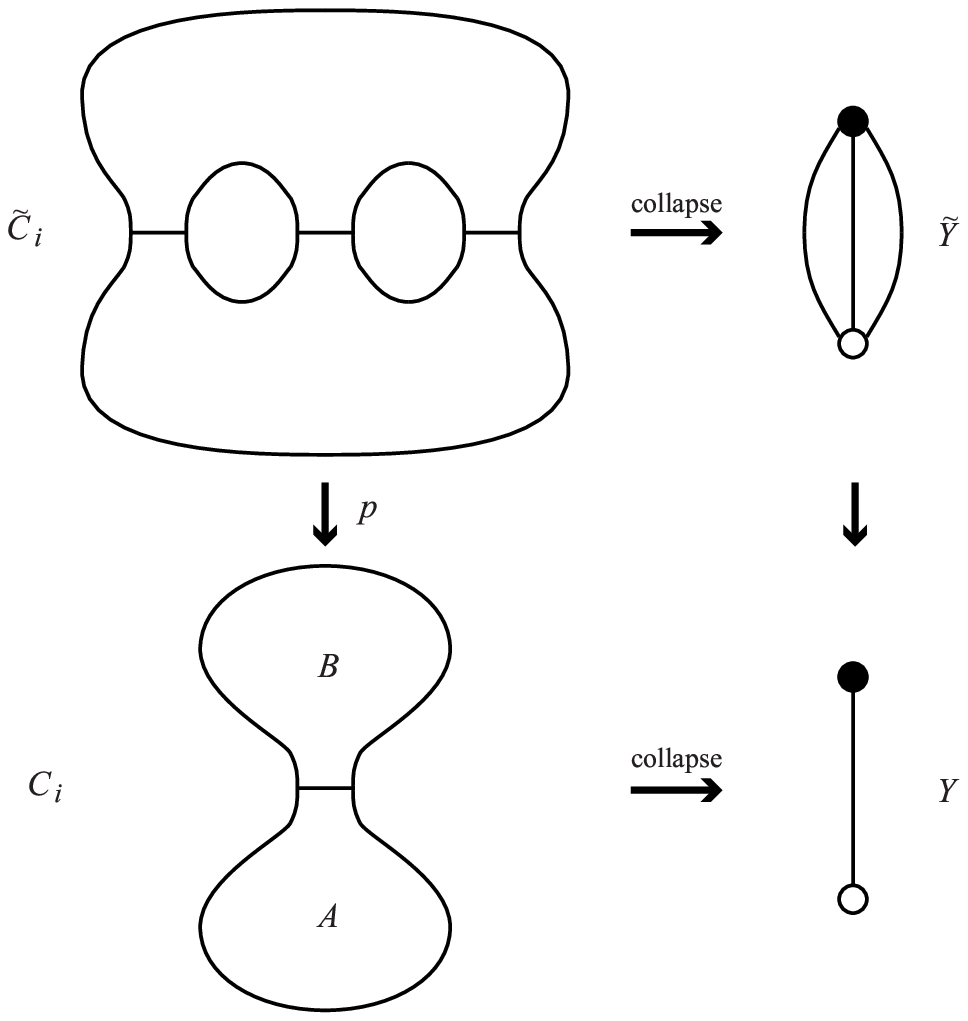,width=3.3in}}
\vskip 18pt
\centerline{Figure 6.}

\noindent {\sl Claim.} For any vertex in $\tilde P$ that maps to $u$,
the number of edges of $\tilde P$ to which
it is incident is
$${[i_\ast \pi_1U: i_\ast \pi_1 U \cap K_i] \over 
[i_\ast \pi_1E: i_\ast \pi_1 E \cap K_i]}.$$

Since $p \colon \tilde C_i \rightarrow C_i$ is a
regular cover, it has a group of covering transformations.
This descends to a group action on $\tilde Y$ which
acts transitively on the vertices that map
to $u$ and acts transitively on the edges that
map to $e$. Hence, we need only show that, in
$\tilde Y$, the number of $e$-labelled edges
and the number of $u$-labelled vertices are
in the ratio given by the above formula. In
fact, we will show that the number of these
edges and vertices are, respectively, 
$${ [J_i:K_i] \over [i_\ast \pi_1E: i_\ast \pi_1 E \cap K_i] },
\qquad
{ [J_i:K_i] \over [i_\ast \pi_1U: i_\ast \pi_1 U \cap K_i]},$$
which will prove the claim. Consider one such edge.
It corresponds to a component $Z$ of $p^{-1}(E)$.
Now, $p|_Z \colon Z \rightarrow E$ is a covering map.
It corresponds to the subgroup of $\pi_1(E)$ 
consisting of based loops in $E$ that lift to
loops in $\tilde C_i$. This subgroup is
$i_\ast^{-1}K_i$, and hence the
degree of $p|_Z$ is $[\pi_1 E: i_\ast^{-1}K_i] = [i_\ast \pi_1E: i_\ast \pi_1 E \cap K_i]$.
But the degree of $p$ is $[J_i:K_i]$, and so
the number of covering translates of $Z$ in $\tilde C_i$, which is
the number of $e$-labelled edges, is given by the
required formula. The same argument applies to the
$u$-labelled vertices, proving the claim.

Now, $$d(i_\ast \pi_1U/ i_\ast \pi_1 U \cap K_i) -
d(i_\ast \pi_1E / i_\ast \pi_1 E \cap K_i) 
=wt(u) - wt(e) \geq g(u) > 1.$$
Hence $i_\ast \pi_1E / i_\ast\pi_1 E \cap K_i$ has
index more than two in $i_\ast \pi_1U/ i_\ast \pi_1 U \cap K_i$.
Thus, the valence of any $u$-labelled vertex in $\tilde P$
is more than two. The same argument
applies to $v$. Hence, every edge of $\tilde P$ has
valence at least two, and at least two vertices
have valence at least three. Therefore,
$\chi(\tilde P) < 0$, and so $\chi(\tilde Y) < 0$.
Thus, $\pi_1(\tilde Y)$ is a free non-abelian group.
But the collapsing map $\tilde C_i \rightarrow \tilde Y$
induces a surjective homomorphism
$K_i \rightarrow \pi_1(\tilde Y)$.  
This proves Theorem 4.1 and hence Theorem 1.2. 
$\square$

Theorem 1.2 can be strengthened to give the following,
rather technical result. Unlike Theorem 1.2, this makes
a hypothesis about just one triple $H \geq J \geq K$
of finite index normal subgroups of $G$.

\noindent {\bf Theorem 4.2.} {\sl Let $G$ be a group
with finite presentation $\langle S | R \rangle$.
Let $L$ be the sum of the lengths of the relations $R$,
which we assume is at least one. Suppose that
$G$ has finite index normal subgroups $H \geq J \geq K$
such that $H/J$ is abelian and
$$d(J/K) > \max \left\{ {48L|S| \over \lfloor([H:J]-1)^{[G:H]^{-1}|S|^{-1}}\rfloor}
+ {16L |S| \over [H:J]}, 16|S| + 8 L^2 \right\}.$$
Then $K$ admits a surjective homomorphism onto a 
non-abelian free group.}

\noindent {\sl Proof.}
All that one needed to make the proof of Theorem 4.1
work was the existence of a positive real number
$\ell$ such that inequalities (1), (2) and (3)
all hold. If we set $\ell$ to be $d(J/K)/[G:J]$,
then (1) automatically is satisfied, and so
we have reduced this to two inequalities.
Inequality (3) becomes $d(J/K) \geq 16|S| + 8 L^2$,
which is part of our hypothesis.
To guarantee that inequality (2) holds, we
use the inequalities in the proof of Theorem 3.2
(where $\Sigma$ is as defined in Lemma 3.3):
$$\eqalign{
{w(X(G/J;S)) \over [G:J]} 
&\leq
{w(X(H/J; \Sigma)) \over [G:J]} +
{2 |S| [G:H] \over [G:J]}\cr
&\leq {6|\Sigma| \over [G:H]\lfloor([H:J]-1)^{|\Sigma|^{-1}}\rfloor}
+ {2 |S| \over [H:J]} \cr
&\leq {6|S| [G:H] \over [G:H]\lfloor([H:J]-1)^{[G:H]^{-1}|S|^{-1}}\rfloor}
+ {2 |S| \over [H:J]} \cr
&< {d(J/K) \over 8L},}$$
by our hypothesis. Hence, the proof of Theorem 4.1
gives a surjective homomorphism from $K$ onto
a free non-abelian group. $\square$

\vfill\eject
\centerline{\caps 5. Groups with more generators
than relations}
\vskip 6pt

In this section, we prove Corollaries 1.3 and 1.4.

\noindent {\sl Proof of Corollary 1.3.}
Let $\langle X | R \rangle$ be a presentation
for $G$ with $|X| - |R| > 1$. We
define a nested sequence of 
finite index normal subgroups
$G = G_1 \geq G_2 \geq \dots $ recursively,
by setting $G_{i+1} = [G_i,G_i](G_i)^i$.
Note that each is a characteristic
subgroup of its predecessor, and so
they are all normal in $G$.

The Reidermeister-Schreier process provides
a presentation for $G_i$ with a total of $(|X| - 1)[G:G_i] + 1$
generators and $|R|[G:G_i]$ relations.
Hence, the abelianisation of $G_i$
contains ${\Bbb Z}^{c_i}$ as a summand,
where $c_i = (|X| - 1 - |R|)[G:G_i] + 1$.
Therefore, $G_i/G_{i+1}$ is an abelian
group with $({\Bbb Z}/i {\Bbb Z})^{c_i}$ as a summand.
So, (i) of Theorem 1.1 holds.
We now verify (ii) and (iii):
$$d(G_i/G_{i+1}) \geq c_i  = (|X| - 1 - |R|)[G:G_i] + 1,$$
which gives (iii), and
$$\log [G_i:G_{i+1}] \geq
\log \left( i^{c_i} \right )
= ((|X| - 1 - |R|)[G:G_i] + 1) \log i,$$
which implies (ii). So, Theorem 1.1
now gives the corollary. $\square$

\noindent {\sl Proof of Corollary 1.4.}
This is similar to the above proof.
Again, let $\langle X | R \rangle$ be
the finite presentation of $G$, with
$|X| - |R| = 1$. Let $w^q$ be the
relation in $R$ that is a proper power.
We may assume that $q$ is prime. Define 
a nested sequence of finite index normal subgroups
$G = G_1 \geq G_2 \geq \dots $ recursively,
by setting $G_{i+1} = [G_i,G_i](G_i)^{p_i}$,
where $p_i$ is the $i^{\rm th}$ prime bigger than $q$.
Again, the Reidermeister-Schreier process gives
a presentation for $G_i$ with $(|X| - 1)[G:G_i] + 1$
generators and $|R|[G:G_i]$ relations.
This has one more generator that relation.
A presentation for $G_i/[G_i,G_i]$ is
obtained by taking the $(|X| - 1)[G:G_i] + 1$
generators, then abelianising,
and then quotienting by the $|R|[G:G_i]$ relations.
Let $H_i$ be the abelian group obtained
by the above procedure, but without quotienting
by the relations that are lifts of $w^q$.
This has ${\Bbb Z}^{[G:G_i]+1}$
as a summand.

We claim that $w \in G_i$ for each $i$.
We will show this by induction. It clearly
holds for $i = 1$. Suppose therefore that
$w \in G_i$. Now, the image of
$w$ in $H_i$ lies in the subgroup of $H_i$ generated by 
$w^q$ and $w^{p_i}$, since $q$ is coprime to
$p_i$. So, $w \in G_{i+1}$, proving the
claim.

Suppose now that, for some $i$,
$w$ has infinite order in $H_i$.
Let $n$ be the largest integer (possibly zero) such that
$w \in q^n H_i$. Set $K$ to be $[G_i, G_i] (G_i)^{q^{n+1}}$.
Then $G/K$ is a finite homomorphic image of $G$
in which $w$ is non-trivial.
A lemma of Baumslag and Pride in [2] then
asserts that $K$, and hence $G$, is large.
(Alternatively, one can avoid
using this lemma of Baumslag and Pride, by showing
directly that $K$ has a presentation
with at least two more generators than relations,
and then using Corollary 1.3.)

Hence, we may assume that $w$ has finite order in $H_i$
for each $i$. So, when we quotient $H_i$ by the relations that
are lifts of $w^q$, a ${\Bbb Z}^{[G:G_i]+1}$ summand remains. 
Therefore, $G_i/G_{i+1}$ contains $({\Bbb Z}/p_i{\Bbb Z})^{[G:G_i]+1}$
as a summand. So,
$$d(G_i/G_{i+1}) \geq [G:G_i] + 1,$$
$$\log [G_i:G_{i+1}] \geq
\log \left( {p_i}^{[G:G_i]+1} \right )
= ([G:G_i]+1) \log p_i.$$
Applying Theorem 1.1 then gives the corollary. $\square$

\vskip 18pt
\centerline{\caps 6. Groups with positive virtual first Betti number}
\vskip 6pt

In this final section, we establish Theorem 1.5.
The implication $(1) \Rightarrow (2)$ is straightforward.
For, if some finite index normal subgroup $G_1$ of $G$ has infinite
abelianisation, we may then find finite index
subgroups $G_1 \geq G_2 \geq \dots$, each normal
in $G_1$, so that $[G:G_i]$ grows as fast
as we like. In particular, we may ensure
that $[G_i:G_{i+1}]$ grows super-exponentially
as a function of $[G:G_i]$. 

The implication $(2) \Rightarrow (1)$ 
is a consequence of the following proposition.

\noindent {\bf Proposition 6.1.} {\sl Let $G$
be a finitely presented group. Then, there is
a constant $k$ with the following property.
If $H$ is a finite index subgroup of $G$,
then either its abelianisation $H/H'$ is infinite or
$|H/H'| \leq k^{[G:H]}$.}

\noindent {\sl Proof.} Pick some finite presentation
for $G$. Let $c$ denote the maximal length
of any of its relations, and set $d$ to be the
number of its generators. Let $K$ be the 2-complex
arising from this presentation. Let $\tilde K \rightarrow K$
be the cover corresponding to a finite index subgroup $H$.
For $j = 0,1,2$, let $C_j$ be the chain
group generated by the $j$-cells of $\tilde K$,
and let $\partial_2 \colon C_2 \rightarrow C_1$
and $\partial_1 \colon C_1 \rightarrow C_0$ be
the boundary maps. Then $H/H'$ equals
${\rm Ker}(\partial_1) / {\rm Im}(\partial_2)$.
Pick a maximal tree in the 1-skeleton of
$\tilde K$. From this, one can construct a
basis of ${\rm Ker}(\partial_1)$, where each
element corresponds to an edge not in the
tree, and is a loop formed by the edge
and the path in the tree joining its
endpoints. Thus, ${\rm Ker}(\partial_1)$
is a product of $[G:H](d-1) + 1$ copies
of ${\Bbb Z}$. Now, ${\rm Im}(\partial_2)$
is spanned by the images of the 2-cells.
Each runs over at most $c$ 1-cells in
$\tilde K$, and so maps to a vector
in $C_1$ with length at most $c$. We may
suppose that $|H/H'|$ is finite, and hence
that ${\rm Im}(\partial_2)$ is a finite
index subgroup of ${\rm Ker}(\partial_1)$. Pick a minimal
collection ${\cal C}$ of 2-cells such
that $\partial_2(\langle {\cal C} \rangle)$
has finite index in ${\rm Ker}(\partial_1)$. Then, the
images of these 2-cells span a paralleliped
in ${\rm Ker}(\partial_1)$. Denote its volume by $V$, which equals
$|{\rm Ker}(\partial_1) / \partial_2(\langle {\cal C} \rangle)|$.
But $|H/H'|$ is a quotient
of ${\rm Ker}(\partial_1)/ \partial_2(\langle {\cal C} \rangle)$,
and so $|H/H'|$ is at most $V$. Since each 2-cell
in ${\cal C}$ maps to a vector in ${\rm Ker}(\partial_1)$ with length
at most $c$, we deduce that the volume
$V$ is at most $c^{[G:H](d-1) + 1} \leq c^{d [G:H]}$. Thus,
the proposition is proved by setting $k$
to $c^d$. $\square$

\vskip 18pt
\centerline{\caps References}
\vskip 6pt

\item{1.} {\caps B. Baumslag, S. Pride,} 
{\sl Groups with two more generators than relators,}
J. London Math. Soc. (2) 17 (1978) 425--426.

\item{2.} {\caps B. Baumslag, S. Pride,} 
{\sl Groups with one more generator than relators,}
Math. Z. 167 (1979) 279--281.

\item{3.} {\caps D. Gabai,} {\sl Foliations
and the topology of 3-manifolds, III},
J. Differential Geom. {\bf 26} (1987) 479--536.

\item{4.} {\caps M. Gromov}, {\sl Volume and
bounded cohomology}, I.H.E.S. Publ. Math. 56 (1982) 5--99.

\item{5.} {\caps M. Lackenby}, {\sl Expanders, rank and
graphs of groups}, To appear in Israel J. Math.

\item{6.} {\caps A. Lubotzky}, {\sl Discrete Groups,
Expanding Graphs and Invariant Measures}, Progr.
in Math. 125 (1994)

\item{7.} {\caps A. Lubotzky, B. Weiss}, {\sl Groups and expanders},
Expanding graphs (Princeton, 1992) 95--109, DIMACS Ser.
Discrete Math. Theoret. Comput. Sci, 10, Amer. Math.
Soc., Providence, RI, 1993.

\item{8.} {\caps A. Lubotzky, R. Zimmer}, {\sl
Variants of Kazhdan's property for subgroups of semisimple groups}, 
Israel J. Math. {\bf 66} (1989) 289--299.

\item{9.} {\caps M. Scharlemann}, {\sl Thin position in the 
theory of classical knots,} To appear in the Handbook
of Knot Theory, Elsevier.

\item{10.} {\caps R. St\"ohr}, {\sl Groups with one more
generators than relators}, Math. Z. 182 (1983) 45--47.

\vskip 12pt
\+ Mathematical Institute, Oxford University, \cr
\+ 24-29 St Giles', Oxford OX1 3LB, UK. \cr

\end